\documentclass[twocolumn,letterpaper]{IEEEtran}
\usepackage[letterpaper,top=.9in,bottom=0.75in,right=0.75in,left=0.75in]{geometry}
\label{gpdyn}

\usepackage{psfrag}
\usepackage{amsmath}
\usepackage{amssymb}
\usepackage{graphicx}
\usepackage{bm} % this is for having mathematical things
\usepackage{bbm} % some math symbols
\usepackage{subfigure}
\usepackage{verbatim}
\usepackage[thinlines,thiklines]{easybmat}
\usepackage{latexsym}
\usepackage[dvipsnames]{xcolor}
\usepackage{cite}
\usepackage{stmaryrd}
\usepackage{multirow}
\usepackage{textcomp}
\usepackage{cleveref} % this is for autumatically typing labels, \cref and \Cref
\crefname{subsection}{subsection}{subsections}
\Crefname{appsec}{appendix}{appendices}
\usepackage{tablefootnote} % for footenotes inside tables
\usepackage{multirow} % for multirow
\usepackage{url} % to show the url
\usepackage{float}
\usepackage{authblk} % for author info
\usepackage[ruled]{algorithm2e}
\newbox\tempbox

\def\cov{\mathop{\rm cov}\nolimits}

\def\diag{\mathop{\rm diag}\nolimits}

\begin{document}

\bibliographystyle{IEEEtran}

\title{\bf \large Gaussian Process (GP)-based Learning Control of Selective Laser Melting Process~\thanks{The work was supported in part by the US National Science Foundation under award CNS-1932370.}}
	
\author{Farshid Asadi~\thanks{F. Asadi, A. Olleak, J. Yi, and Y. Guo are with the the Department of Mechanical and Aerospace Engineering, Rutgers University, Piscataway, NJ 08854 USA (e-mail: fa416@scarletmail.rutgers.edu).}\hspace{-4pt}, Alaa Olleak, Jingang Yi, and Yuebin Guo}
	
\maketitle

\begin{abstract}
Selective laser melting (SLM) is one of emerging processes for effective metal additive manufacturing. Due to complex heat exchange and material phase changes, it is challenging to accurately model the SLM dynamics and design robust control of SLM process. In this paper, we first present a data-driven Gaussian process based dynamic model for SLM process and then design a model predictive control to regulate the melt pool size. Physical and process constraints are considered in the controller design. The learning model and control design are tested and validated with high-fidelity finite element simulation. The comparison results with other control design demonstrate the efficacy of the control design.    
\end{abstract}

%\vspace{-0.5mm}
\section{Introduction}

Metal additive manufacturing (AM) builds three-dimensional parts by melting metal powders layer-by-layer using laser or other heating sources. Selective laser melting (SLM) is one of the popular metal AM processes. Real-time feedback control of metal AM is an enabling tool for reliable, robust, high-quality process~\cite{TapiaJMSE2014}. Although repetitive feature is used for control of metal AM~\cite{TapiaJMSE2014}, few work exist for real-time control of SLM due to its complex process features and lack of enabling, effective and reliable sensing techniques. The goal of this paper is to develop a data-driven SLM process modeling and control design. 

Control of metal AM has been reported in past decade. For example, built on repetitive nature of the AM process, iterative learning control is an effective way to improve the fabrication of the metal AM processes (e.g.,~\cite{SammonsCST2019,ShkorutaCASE2019}). The dynamic models for directed energy deposition (DED) processes are mainly built on mass conversation and energy balance of the melt pool (e.g.,~\cite{WangJMSE2017,SammonsJDSMC2019}). These physics-based model become complex for modeling SLM process because of dynamic behavior of the powder melt pool, non-uniform heat dissipation from the melt pool to the surrounding materials, and various process parameters and scanning patterns. The recent work in~\cite{WangAM2020} built a simplistic physics-based model for real-time track-by-track control of SLM. The effect of previous scanning track is considered in the next track control. Several important assumptions are made for the melt pool geometry to derive the model using the energy balance principle. Finite element (FE) method~(e.g.,\cite{OlleakJMD2020}) and multi-scale multi-physics models (e.g.,~\cite{KingMST2015,YanCM2018}) are used to obtain high-fidelity simulation for SLM process. It is however, impossible to use these models to design real-time process control due to their high computational costs. The recent work in~\cite{Wang2020}, used data obtained from FE analysis to derive a linear model of the process. In their work, the authors used a repetitive control algorithm to control melt pool width for a multi-track case in simulated FE software.

Surrogate models are developed to optimize the SLM process parameters for stable melt pool and process windows of the laser power and velocity maps~\cite{TapiaAM2016,TapiaIJAMT2018,YeungAM2019,YeungAM2020}. SLM powder material characterization and process properties (e.g., melt pool geometry, porosity, etc.) are also reported by using the Gaussian process (GP) models~\cite{TapiaAM2016,TapiaIJAMT2018} and it is thus attractive to integrate these development into the real-time control design. However, these data-driven models do not incorporate temporal information and therefore cannot directly be applied to capture dynamics of the SLM process. In this paper, we extend the GP models to capture the dynamics feature of the SLM process and propose a learning-based real-time process control. 

\begin{figure*}[htb!]
	\centering
	\subfigure[]{
		\label{fig_slm_setup}
		\includegraphics[width=.3\textwidth]{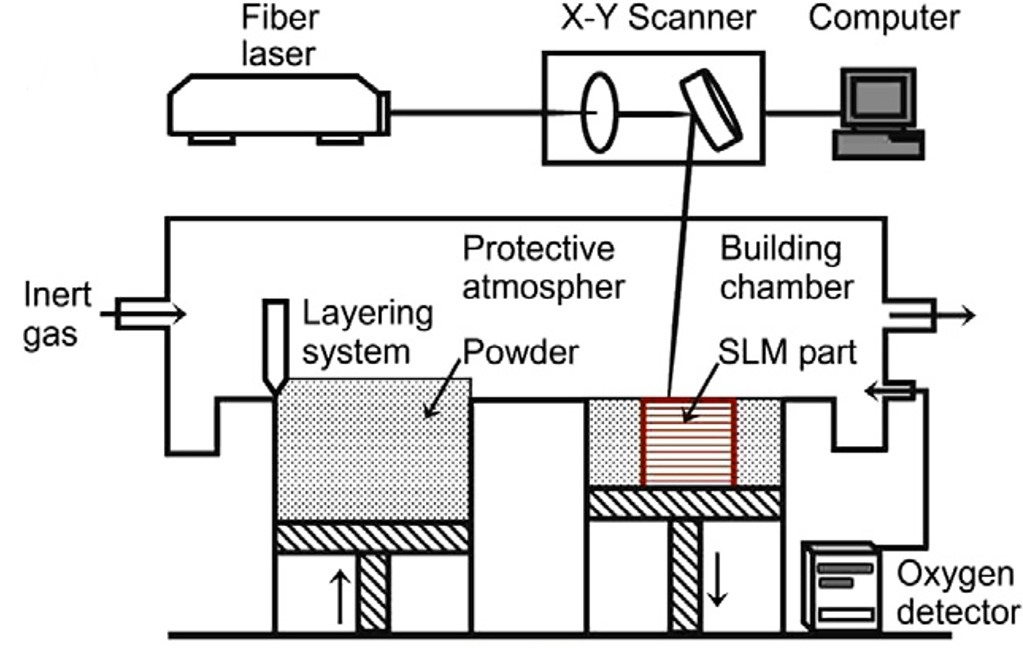}}
	\hspace{0mm}
	\subfigure[]{
		\label{fig_slm_scan_a}
		\includegraphics[width=.3\textwidth]{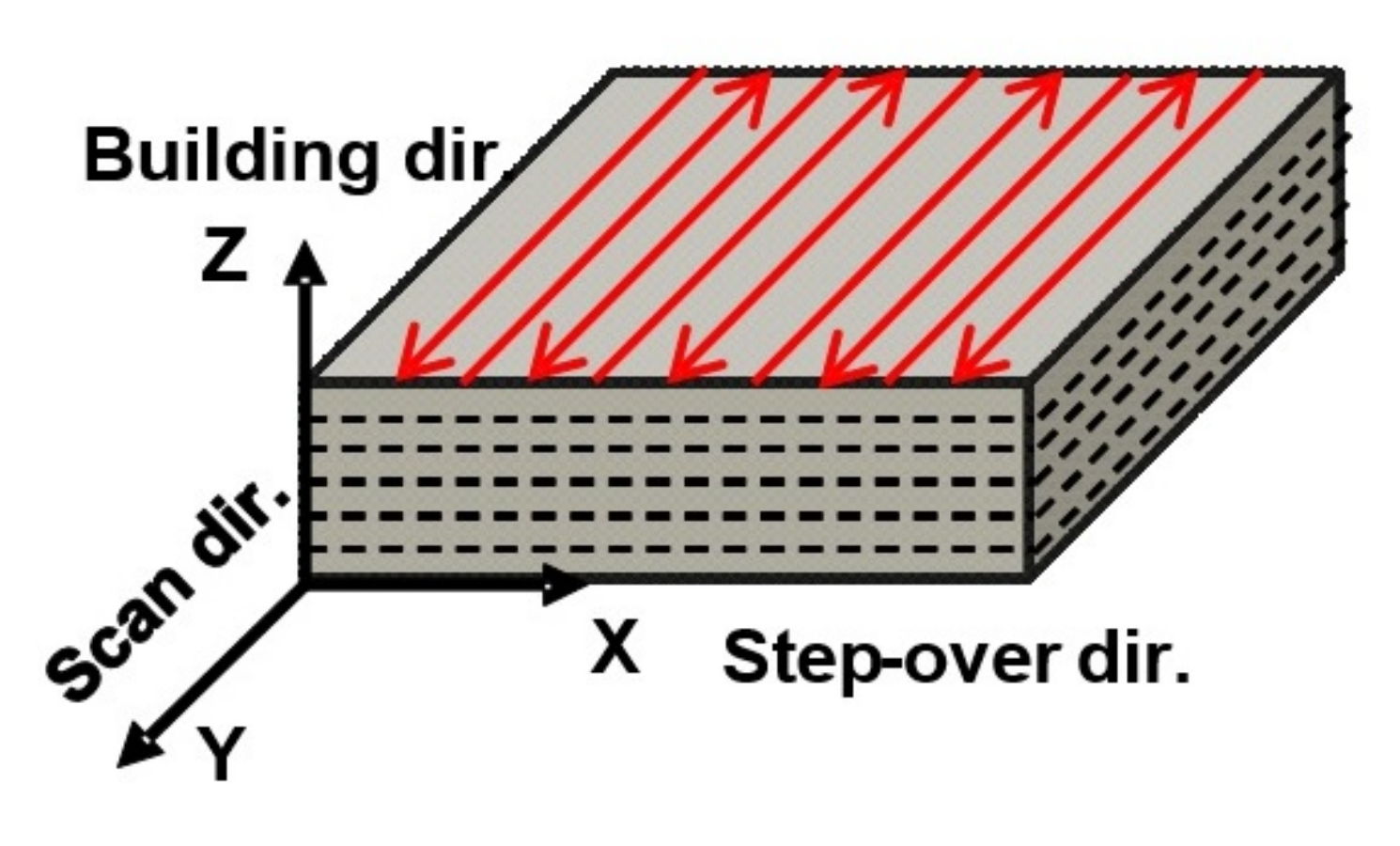}}
	\hspace{0mm}
	\subfigure[]{
		\label{fig_slm_scan_b}
		\includegraphics[width=.3\textwidth]{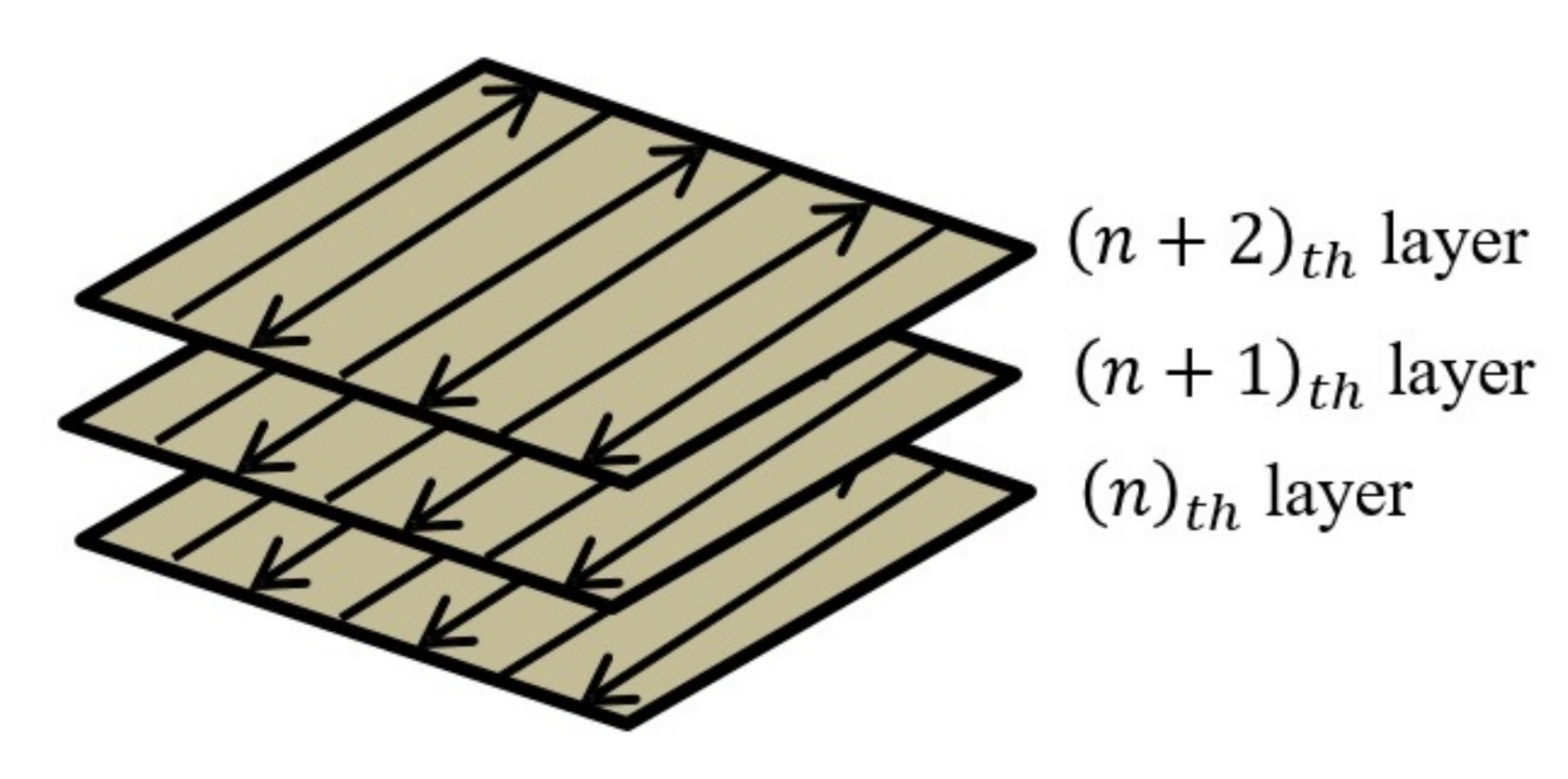}}
	\caption{(a) SLM process experimental setup. (b) Rectangular SLM scanning strategy. (c) Layer-by-layer configuration for SLM process.}
	\label{fig_slm_scan}
\end{figure*}
GP dynamics model (GPDM) and model predictive control (MPC) was proposed and used in various applications~\cite{WangPAMI2008,ChenICRA2019,HanRAL2021}. In this paper, we take advantages of the GP learning model of the SLM dynamics and develop an optimization-based process control. The data-driven GP model is built on the calibrated FE simulation with limited experimental data. The main contribution of the work lies in the data-driven GP model and the design of real-time MPC control of the stable melt pool in SLM.

\section{Data-Driven SLM Process Models}
\label{model}
\subsection{SLM Process Dynamics}
\label{slmdyn}
SLM is a powder bed fusion (PBF) process where the powder particles absorb the energy delivered by the moving laser beam, start to melt, and form a melt pool. Fig.~\ref{fig_slm_setup} illustrates the experimental setup of a typical SLM process. The laser beam movement is controlled by mirrors that direct the beam in the $x$ and $y$ directions at high speeds and the laser beam is shot on powder particles layers. During a SLM process, a thin layer of materials powder is distributed over the solid substrate by the recoater. Under the laser beam scanning, a desired pattern is fused to the substrate. For perfect bonding, the melt pool must be deeper than the powder layer thickness and wider than hatch spacing (i.e., the distance between center lines of two adjacent tracks). After the entire layer is scanned according to the scanning strategy, the next powder layer is added and the process is repeated until the entire part is built. Fig.~\ref{fig_slm_scan_a} shows one layer scanning SLM process and Fig.~\ref{fig_slm_scan_b} illustrates the multiple layer building structure.

It was reported that over one hundred parameters affect the quality of parts manufactured by SLM, among which the most important parameters are laser power and scanning speed. Furthermore, residual heat effect due to the subsequent scanning tracks might result in high melt pool size variations over the entire layer. Thus, real-time process monitoring and feedback control are needed to achieve a stable and consistent melt pool sizes. Process monitoring can be conducted through fast thermal imaging by using CCD camera, pyrometer, and light sensitive diodes, e.g.,~\cite{doumanidis2001geometry}. These sensing suites can be used to provide information about the temperature distribution, and melt pool size, and therefore to develop a feedback control enhancing the process behavior. In this work, a single layer with multi bidirectional tracks are considered for SLM modeling and control.

\subsection{Gaussian Processes Modeling}
\label{gp}

GP modeling is used to obtain data-driven model of the SLM process dynamics. We here present a short summary of the GP models and readers can refer to~\cite{rasmussen2017gaussian} for further details. GP is a powerful technique to derive analytical models from relatively small data set. We consider a noisy observation of the process as
\begin{equation}
\label{eq_obsform}
y = f(\bm x) + \eta,
\end{equation}
where $y\in \mathbb{R}$ is a scalar output observation, $\bm x \in \mathbb{R}^n$ is an $n$-dimensional independent variable, $f(\cdot): \mathbb{R}^{n} \rightarrow \mathbb{R}$ is the underlying unknown model, and $\eta \sim \mathcal{N}(0,\sigma^2_n)$ is zero mean Gaussian observation noise with variance $\sigma^2_n$. Furthermore, we assume that $f(\cdot)$ has a Gaussian probability distribution prior. Without loss of generality, $f(\cdot)$ is considered as a Gaussian process with the (zero) mean and covariance function
\begin{equation}
\label{eq_gpf}
f(\bm x) \sim \mathcal{GP}(0, k(x,x')),
\end{equation}
where $k(\bm x,\bm x'): \mathbb{R}^{n} \times \mathbb{R}^{n} \rightarrow \mathbb{R}$ is a positive definite covariance function. In this work, a squared exponential function is used with the following form
\begin{equation}
\label{eq_gpcov}
k(\bm x,\bm x') = \sigma^2_f \exp\begin{pmatrix}-\frac{(\bm x-\bm x')^T \bm L^{-1} (\bm x-\bm x')}{2}\end{pmatrix},
\end{equation}
where $\sigma_f$ and $\bm L = \diag\left(l_1, \dots, l_n\right)$, $l_i >0$, $i=1,\ldots,n$, are hyper parameters to be determined.

Having $m$ observation pairs $(y_i,\bm x_i)$, $i=1,\ldots,m$, and the assumption of Gaussian prior over $f(\cdot)$, we construct the following joint distribution of training data set and prediction pair as
\begin{equation}
\label{eq_gpjoint}
\begin{bmatrix}
\bm Y\\f^*
\end{bmatrix} \sim \mathcal{N}
\begin{pmatrix}
\mathbf{0},
\begin{bmatrix} \bm K(\bm X,\bm X) + \sigma^2_n \bm I_m & \bm K(\bm X,\bm x^*)\\\bm K(\bm x^*,\bm X) &K(\bm x^*,\bm x^*)
\end{bmatrix}
\end{pmatrix},
\end{equation}
where $\bm Y = \begin{bmatrix}y_1, \dots, y_m\end{bmatrix}^T$ and $\bm X = \begin{bmatrix} \bm x_1, \dots, \bm x_m\end{bmatrix}^T$ are training data sets, $\bm x^*$ and $f^*$ are prediction pair, and $\bm K(\cdot,\cdot)$ are corresponding covariance given as follows.
\begin{gather*}
\bm K(\bm X,\bm X)=\begin{bmatrix} k(\bm x_1,\bm x_1)& \dots& k(\bm x_1,\bm x_m)\\ \vdots& \ddots& \vdots\\ k(\bm x_m,\bm x_1)& \dots& k(\bm x_m,\bm x_m) \end{bmatrix} \in \mathbb{R}^{m \times m}, 
\end{gather*}
$\bm K(\bm x^*,\bm X)=[k(\bm x^*,\bm x_1) \, \cdots \, k(\bm x^*,\bm x_m)]=\bm K^T(\bm X,\bm x^*)\in \mathbb{R}^{1\times m}$, 
and $K(\bm x^*,\bm x^*)=k(\bm x^*,\bm x^*)$.

Having the joint distribution~(\ref{eq_gpjoint}), we use Bayesian rule to calculate posterior distribution of $f^*$ as a function of $\bm x^*$, conditional on training data set $\bm Y$ and $\bm X$ as 
\begin{equation}
\label{eq_gpposterior}
f^*|\bm X,\bm Y,\bm x^* \sim \mathcal{N}
\begin{pmatrix}
\bar{f}^*, \cov(f^*)
\end{pmatrix},
\end{equation}
where mean $\bar{f}^* = \bm K(\bm x^*,\bm X)\left(\bm K(\bm X,\bm X) + \sigma^2_n \bm I_m\right)^{-1}\bm Y$ and covariance function $\cov(f^*) =K(x^*,x^*)-\bm K(\bm x^*,\bm X)\left(\bm K(\bm X,\bm X) + \sigma^2_n \bm I_m\right)^{-1}\bm K(\bm X,\bm x^*)$.
We then write the marginal likelihood probability for the training data set as
\begin{equation*}
\label{eq_gpmarginal}
\log p(\bm Y|\bm X,\bm \Theta) = -\frac{1}{2}\bm Y^T \bm K_Y^{-1} \bm Y-\frac{1}{2}\log\left| \bm K_Y\right|-\frac{n}{2}\log(2\pi),
\end{equation*}
where $\bm K_Y = \bm K(\bm X,\bm X) + \sigma^2_n \bm I_m$, $\bm \Theta = \{\sigma_f, \sigma_n, \bm L\}$ is the hyper-parameter set, and $| \bm K_Y |$ denotes the determinant of $\bm K_Y$. Maximizing marginal likelihood with respect to $\bm \Theta$ would generate the appropriate hyper parameters.

\begin{figure*}[htb!]
	\centering
	\subfigure[]{
		\label{fig_slm_overview}
		\includegraphics[width=.4\textwidth]{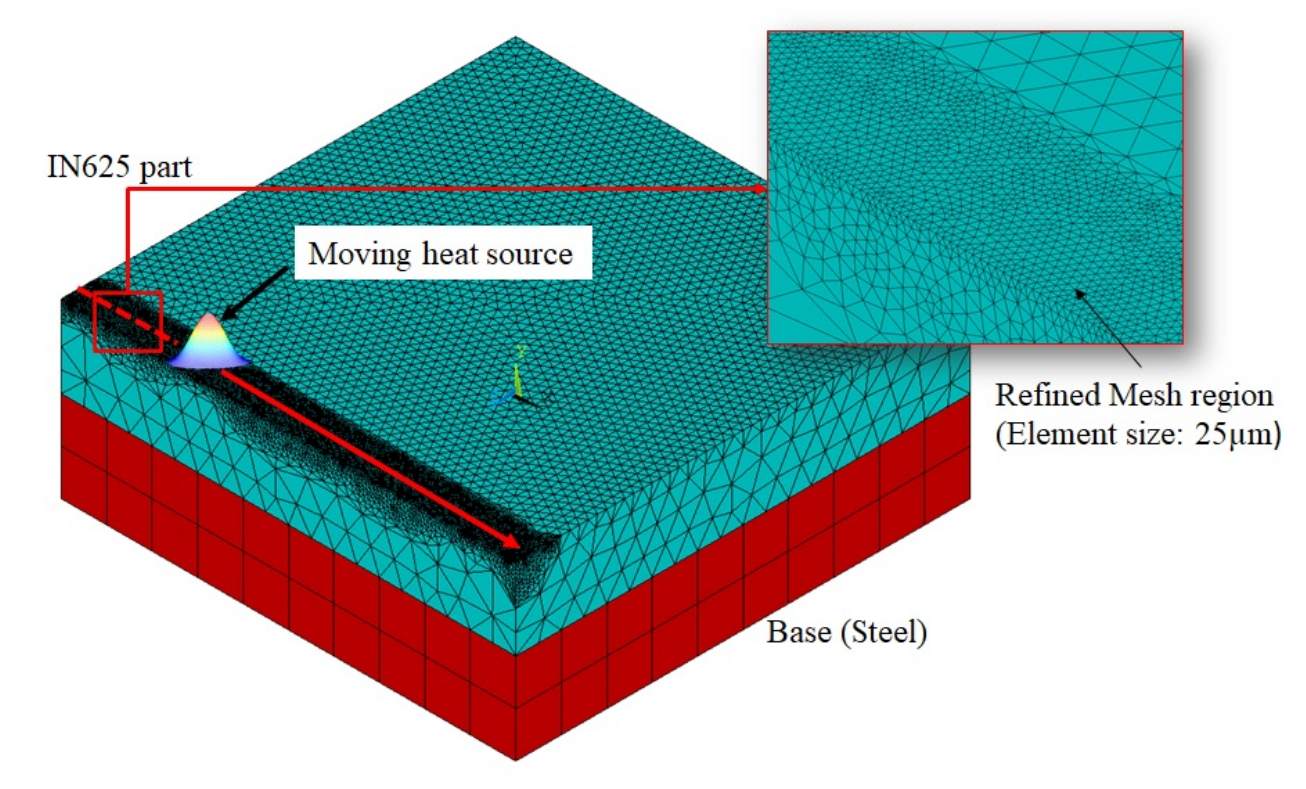}}
	\hspace{10mm}
	\subfigure[]{
		\label{fig_melt_pool_area}
		\includegraphics[width=.4\textwidth]{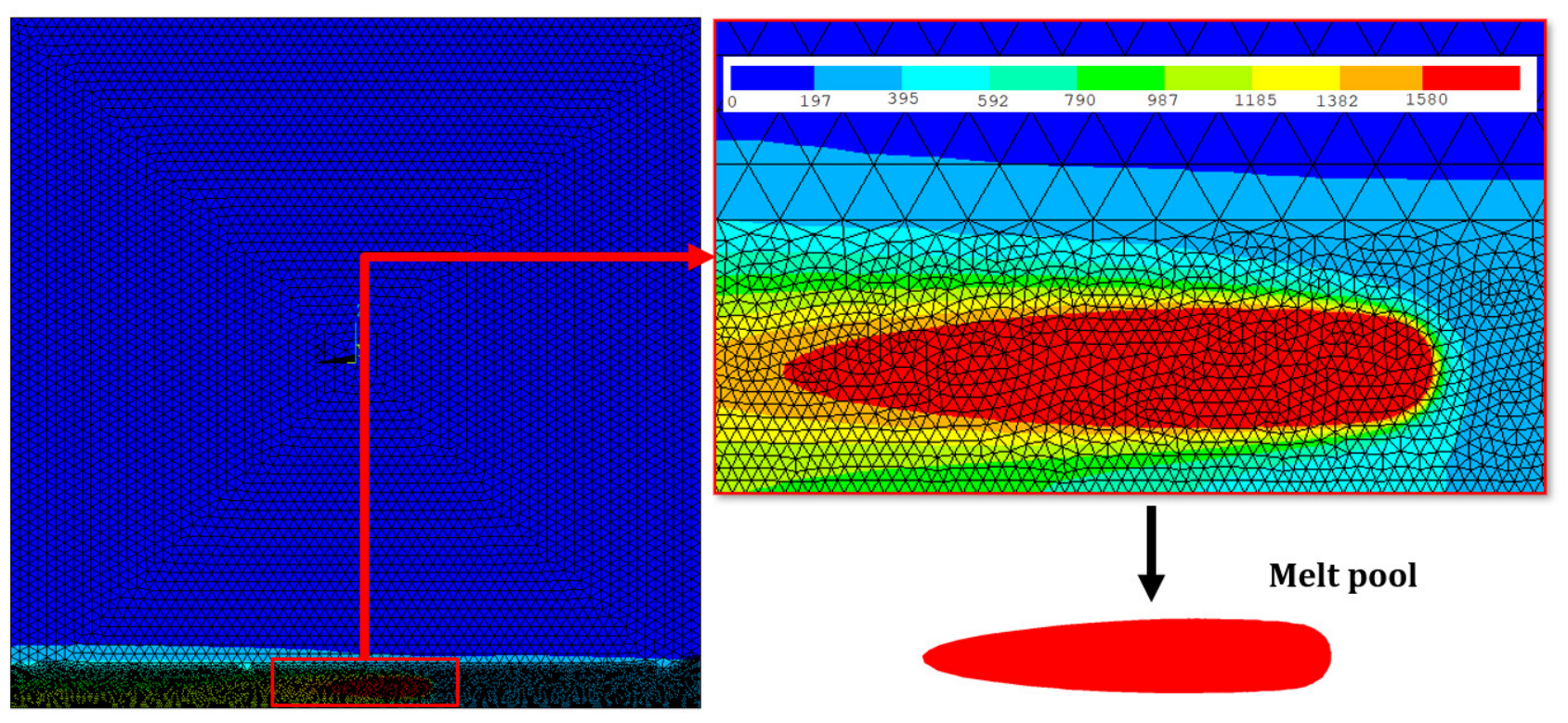}}
	\caption{(a) SLM process FEA overview. (b) Melt pool area extraction from thermal solution.}
\end{figure*}

\subsection{GP-based SLM Process Dynamics}
\label{slmgpdyn}

We consider the GP model for the SLM dynamics. For the SLM process, one of the main performance metrics is the consistent melt pool cross sectional area, denoted as $x(t)$ at time $t$. Among many process parameters, laser power and scanning speed are the most influential process control parameters~\cite{OlleakJMD2020}. Therefore, we consider laser power and scanning speed, denoted as $p(t)$ and $v(t)$ respectively, as control inputs for the SLM dynamics. Additionally, it is observed that substrate initial temperature affects the process states such as cross sectional area~\cite{WangAM2020}. As a result, surface temperature, denoted as $T(t)$ at time $t$, of the selected scanning point before laser beam hits, is considered as a measurable disturbance to the process dynamics~\footnote{At each sampling time, temperature of next laser spot (plus laser beam nominal radius) is used for this purpose.}. This treatment is similar to that in~\cite{WangAM2020}. Therefore, SLM process dynamics is formulated in discrete time as follows
\begin{equation}
\label{eq_dynform}
x_{k+1} = f(x_{k},T_k,p_k,v_k),
\end{equation}
where $x_k$, $T_k$, $p_k$, and $v_k$ are the melt pool cross sectional area, surface temperature of melt pool surrounding, laser power, and laser scanning speed at the $k$th step, respectively, $k \in \mathbb{N}$. 

GP is used to estimate $f$ in~(\ref{eq_dynform}) using training data set. To train the GP models, process data sets are necessary and in many cases, experiments are needed to collect these data sets. However, running SLM experiments is expensive and in many cases, high-fidelity numerical simulation such as finite element analysis has been extensively used as an alternative means to study SLM processes.  

\subsection{Finite Element Analysis}
\label{FE}

Finite Element (FE) analysis is widely used to simulate and predict the thermal history during SLM process. In this work, FE models are developed to provide information about the melt pool size and temperature history to build the GP model and validate the control design. At the $i$th node in the FE model, the addition of heat by the laser beam is represented using the volumetric heat source considering the Gaussian profile for the beam shape, that is, the heat flux is calculated as
\begin{equation}
\label{eq_laser}
Q_{i}=\frac{6\sqrt{3}a p_k}{r_{s}^2c\pi\sqrt{\pi}}e^{-3\left[\bigl(\frac{\Delta x_{i}}{r_{s}}\bigr)^2+\bigl(\frac{\Delta y_{i}}{r_{s}}\bigr)^2+\bigl(\frac{\Delta z_{i}}{c}\bigr)^2\right]},
\end{equation}
where $\Delta x_i$, $\Delta y_i$, and $\Delta z_i$ are the distances between the $i$th node and laser beam center at the $k$th step, $a$ is the absorptivity, $p_k$ is the laser power at the $k$th step, $r_s$ is the laser beam radius, and $c$ is the penetration depth. For computational efficiency and accuracy, the adaptive remeshing framework in~\cite{OlleakJAMT2020} along with tetrahedral mesh are utilized, where a fine mesh of average size $25\,\mu m$ is used at the high temperature gradient regions (e.g., melt pool). Fig.~\ref{fig_slm_overview} illustrates the mesh configuration during scanning the first track. The implicit Ansys\textsuperscript{\textregistered} MAPDL\textsuperscript{\textregistered} solver is used to solve the thermal problem. The movement of the laser beam is simulated by multiple time steps, where the single step represents laser movement distance. Fig.~\ref{fig_melt_pool_area} illustrates the melt pool area and temperature distribution in the FE model. 

\section{SLM Process Control}
\label{control}
Linear time varying MPC is used in the present work for its constraint handling capability and relative low computational cost. Considering the SLM dynamics~(\ref{eq_dynform}), we denote the GP estimate of nonlinear function $f$ as $\bar{f}_{GP}$ and the control input vector $\bm u_k=[p_k \; v_k]^T \in \mathbb{R}^2$. The discrete-time linearized system is then given in the following form,
\begin{equation}
\label{eq_lind}
x_{k+1} = A_d x_k + \bm B_d \bm u_k,
\end{equation}
where $A_d$ and $\bm B_d$ are calculate as
\begin{equation*}
A_d = \frac{\partial \bar{f}_{GP}}{\partial x}\Big\vert_{(x_k,\bm u_k,T_k)}, \; \bm B_d = \frac{\partial \bar{f}_{GP}}{\partial \bm u}\Big\vert_{(x_k,\bm u_k,T_k)},
\end{equation*}
and $T_k$ is the surface temperature of melt pool surrounding at the $k$th step. 

The MPC problem is formulated as a constrained optimization and at the $k$th step, the formulation is given as 
\begin{subequations}
\begin{align}
\min_{\bm U}\quad&\sum_{i = 0}^{H-1}\left(Q x_{k+i}^2 + \bm u_{k+i}^T \bm R \bm u_{k+i}\right) + Q_f x_{k+H}^2, \\
s.t.\quad& x_{i+1} = A_d x_i + \bm B_d \bm u_i, \\
& x_l\leq x_i \leq x_u, \, p_l\leq p_i \leq p_u, \, v_l\leq v_i \leq v_u, \\
& x_l\leq x_{k+H} \leq x_u, \; i=k,\ldots,k+H-1,
\end{align}
\label{eq_mpc}
\end{subequations}
\hspace{-2.5mm} where  $\bm U=\{\bm u_k,\ldots,\bm u_{k+H}\}$ is the control input set and $H$ is the state prediction horizon. In~(\ref{eq_mpc}), $x_l$ ($x_u$), $p_l$ ($p_u$), and $v_l$ ($v_u$) are the lower (upper) bounds for the melt pool cross sectional area, laser power, scanning speed, respectively. MPC control bounds should be chosen based on the optimal parameters of the process to ensure adequate adhesion between melted powder and the substrate and also to avoid keyhole formation~\cite{lo2019optimized}. 

Although the MPC control design is built on the GP model $\bar{f}_{GP}$, the validation of the control performance is through the FE simulation. The use of controller with FE models was previously implemented using COMSOL\textsuperscript{\textregistered} Multiphysics in~\cite{Wang2020}. While COMSOL\textsuperscript{\textregistered} and MATLAB\textsuperscript{\textregistered} can be integrated through LiveLink\textsuperscript{\textregistered}, this is not an option with Ansys\textsuperscript{\textregistered}. In this work, the temperature solution is written by Ansys\textsuperscript{\textregistered} after each step is solved in batch mode. The Python framework, which is integrated and implemented between simulations steps, is responsible for melt pool area calculation and passing the temperature values at location ahead of the melt pool to the MPC controller. The controller updates the laser power and the FE solver proceeds with these updates for the following simulation step. The process is repeated until scanning the entire powder layer is complete. Detailed results will be demonstrated in the next section. 

\section{Simulation Results}
\label{results}
\subsection{GP Model Performance}

We take the IN625 metal SLM process in~\cite{WangAM2020} as an example to demonstrate the GP model training and validation. The IN625 thermal properties are taken from~\cite{WangAM2020} and the value of the process parameters are shown in Table~\ref{table_laser_param}. The implicit Ansys\textsuperscript{\textregistered} APDL\textsuperscript{\textregistered} solver is used to solve the thermal problem. The simulation step time is $\Delta T=50 $ $\mu$s. 

\renewcommand{\arraystretch}{1.3}
\setlength{\tabcolsep}{0.04in}
\begin{table}[tp!]
	\centering
	\caption{Process parameters used in IN625 metal SLM simulation.}
	\label{table_laser_param}
	\begin{tabular}{|c|c|c|c|c|c|c|}
		\hline\hline
		Parameter & $p$ (W) & $v$ (mm/s) & Hatch space (mm) & $r_s$ ($\mu$m) & $a$ & $c$ ($\mu$m) \\ \hline 
		Value & $250$ & $800$ & $0.1$ & $50$ & $0.4$ &	$3$ \\
		\hline\hline
	\end{tabular}
\end{table}

\begin{figure*}[t]
	\centering
	\subfigure[]{
		\label{fig_test_profile}
		\includegraphics[width=.3\textwidth]{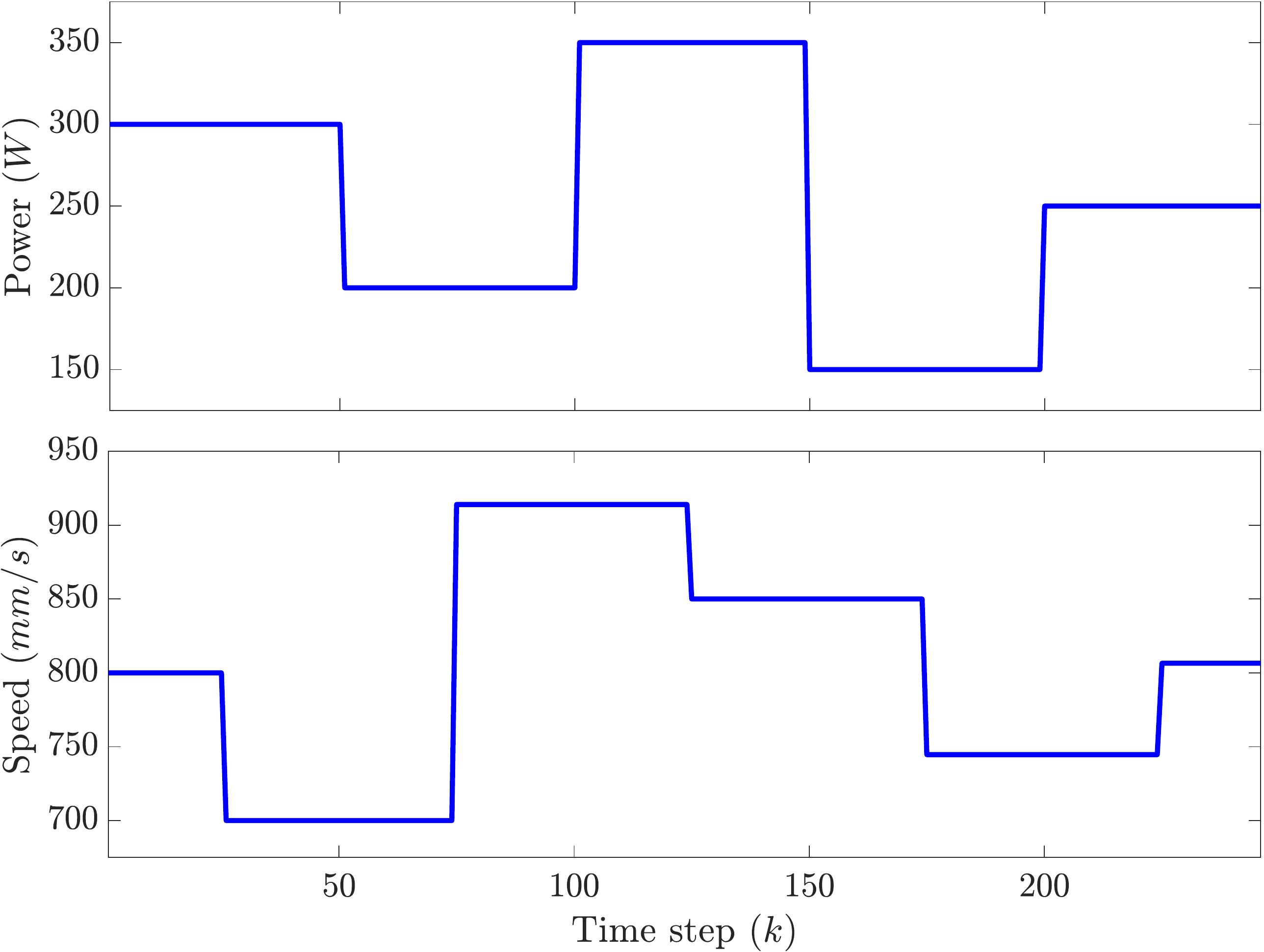}}
	\subfigure[]{
		\label{fig_test_area}
		\includegraphics[width=.3\textwidth]{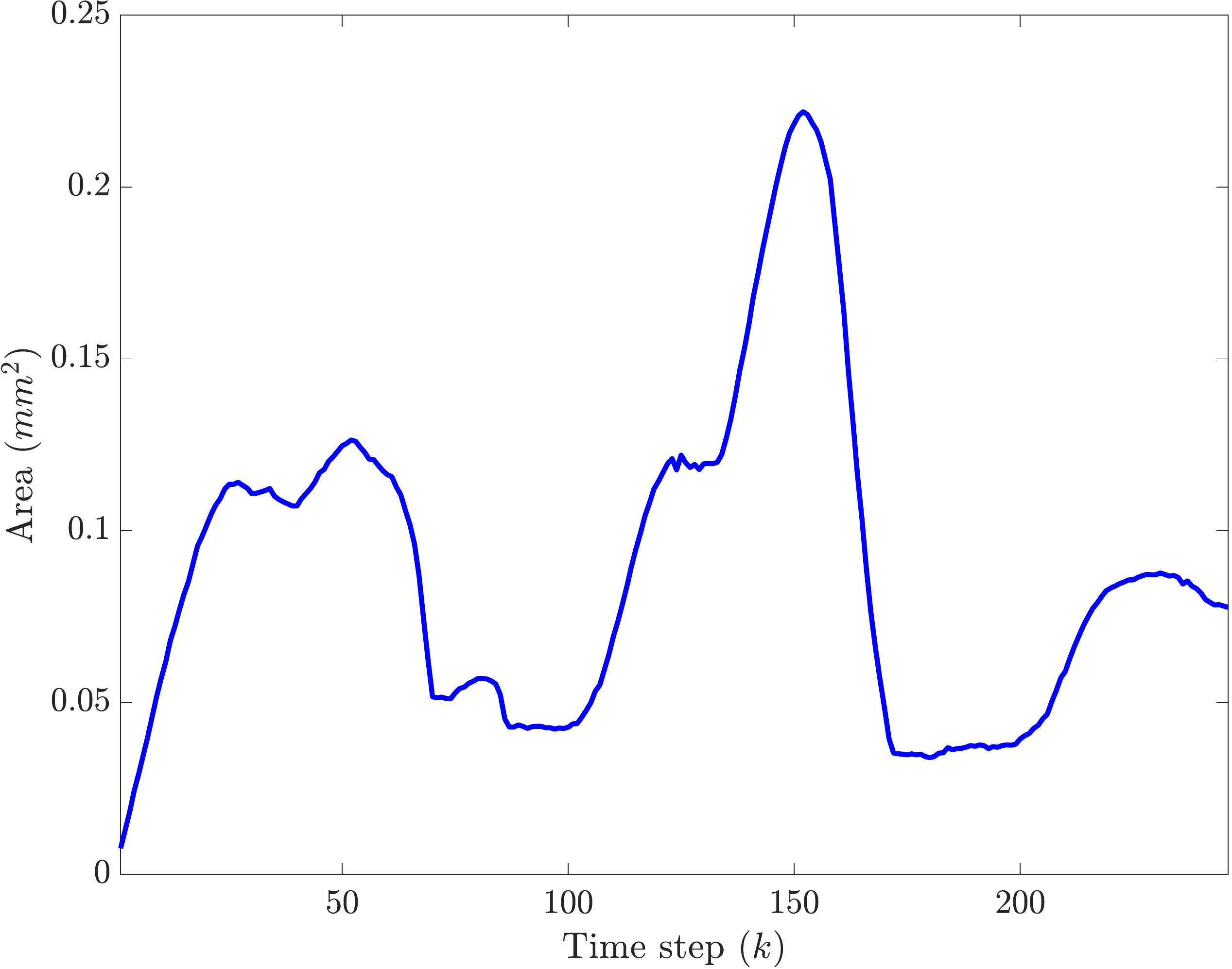}}
	\subfigure[]{
		\label{fig_test_init_temp}
		\includegraphics[width=.3\textwidth]{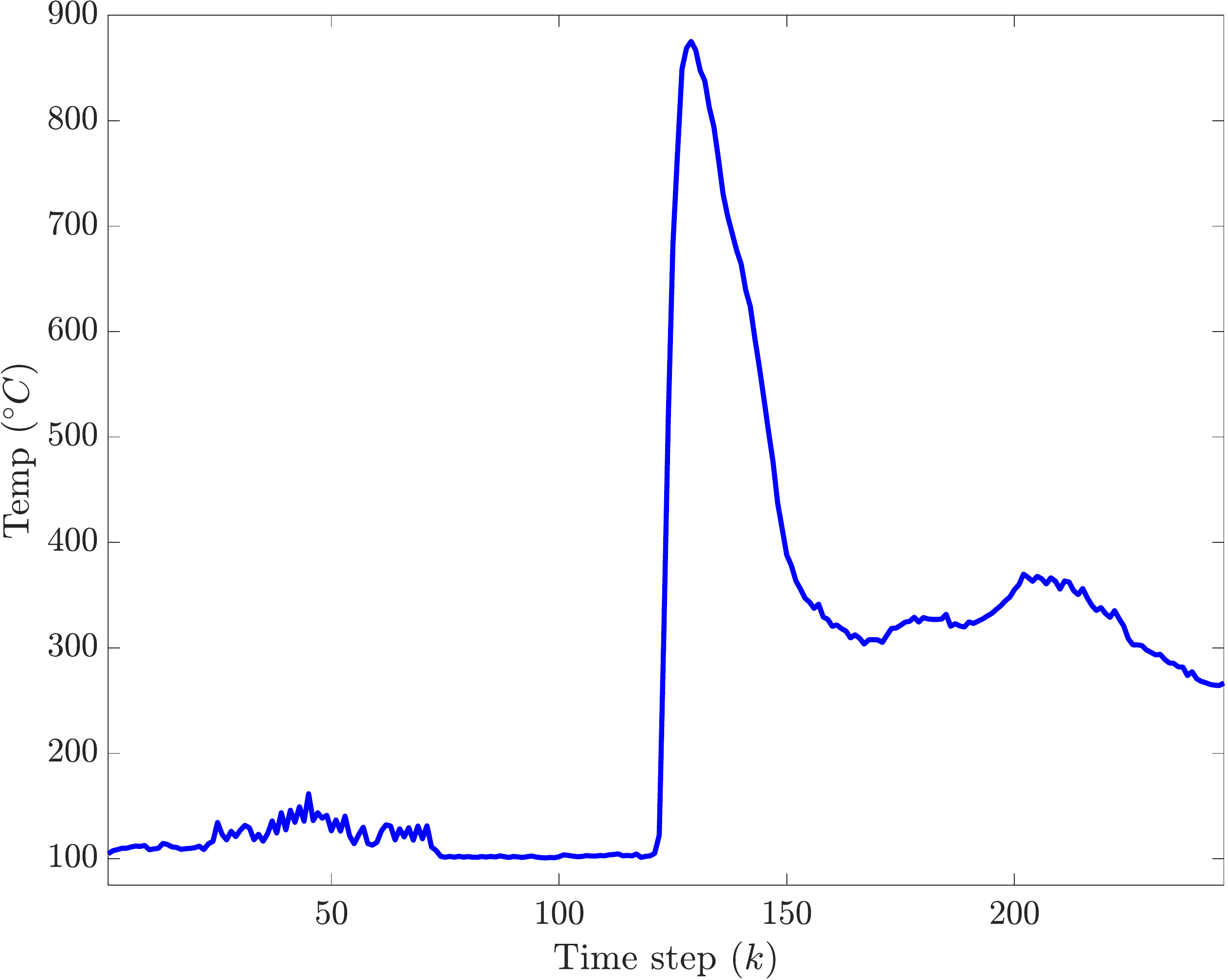}}
	\caption{FE simulation profiles and results for Case \#9 for GP construction. (a) Laser power $p(t)$ and scanning speed $v(t)$ profiles. (b) Melt pool cross sectional area $x(t)$. (c) Surface temperature $T(t)$ of the melt pool surrounding.}
\label{fig_test}
\end{figure*}

\begin{figure*}[htb!]
	\centering
	\subfigure[]{
		\label{fig_training_accuracy}
		\includegraphics[width=.4\textwidth]{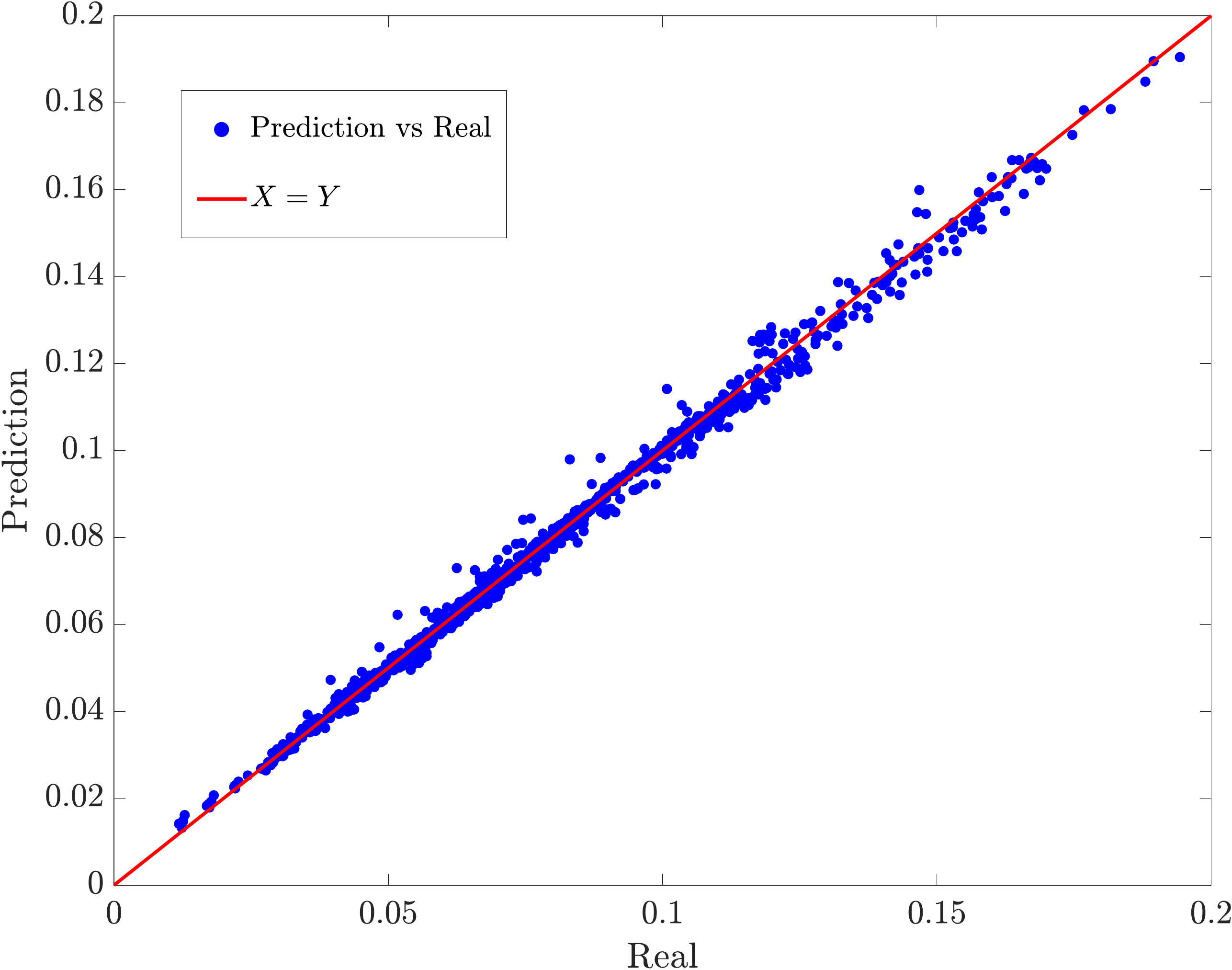}}
	\hspace{10mm}
	\subfigure[]{
		\label{fig_forward_accuracy}
		\includegraphics[width=.4\textwidth]{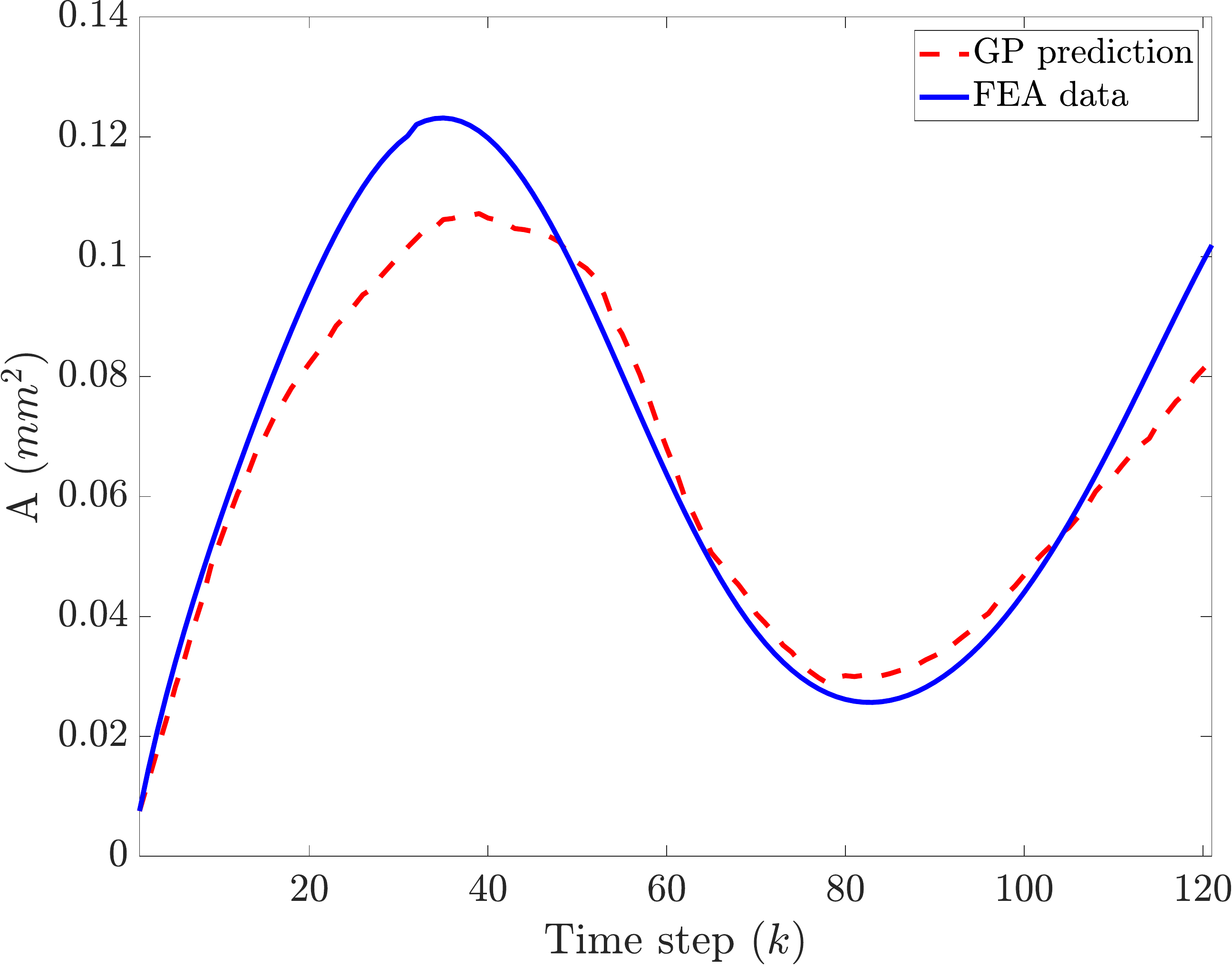}}
	\caption{(a) Comparison results between the trained GP model prediction and the SE simulation using validation data set. (b) Comparison results of the GP-based model prediction over time and the FE simulation of the melt pool cross sectional area for Case \#2 profile.}
\end{figure*}

In order to collect representative data for training the GP model, four single track and five double track simulations are conducted with various combinations of scanning speed and laser power profiles and different hatch spaces (for double track simulations). Table~\ref{table_tests} lists these simulation experiments and their configurations. In all of the test cases, each track length is $5$ mm and the power and scanning speed profiles are constrained within optimal process parameters for SLM process (e.g.,~\cite{lo2019optimized}). The scanning speed and laser power, along with melt pool cross sectional area and scanning point initial temperature, are recorded at sampling time $\Delta T$. As an example, Fig.~\ref{fig_test} shows the laser power $p(t)$ and scanning speed $v(t)$, cross sectional area $x(t)$ and melt pool temperature profiles for Case \#9. 

\begin{table}[hbp!]
	\centering
	\caption{FE simulation summary for GP model training}
	\label{table_tests}
	\begin{tabular}{| l | c | c | c|}
		\hline\hline
		Case No.& Power& Scanning speed & Hatch space (mm) \\ \hline
		Single track&&&\\
		\hline
		\# $1$&Constant & Constant & N/A \\\hline
		\# $2$&Sinusoidal & Constant & N/A \\\hline
		\# $3$&Constant & Sinusoidal & N/A \\\hline
		\# $4$&Sinusoidal & Sinusoidal & N/A \\\hline
		Double track&&&\\
		\hline
		\# $5$&Constant & Constant & $0.1$  \\\hline
		\# $6$&Constant & Constant & $0.15$ \\\hline
		\# $7$&Constant & Constant & $0.05$ \\\hline
		\# $8$&Sinusoidal & Profile & $0.1$ \\\hline
		\# $9$&Profile & Profile & $0.1$ \\
		\hline\hline
	\end{tabular}
\end{table}

From all nine cases FE simulation, a total of 1649 data sets were collected, among which 100 data sets were used for the GP training and the rest of 1549 data points are used for model validation. Fig.~\ref{fig_training_accuracy} shows the comparison of the GP model prediction and the FE simulation results. The average of the prediction error of the GP model is around $2.06$ \% for the melt pool cross sectional area. Furthermore, coefficient of determination $R^2 = 99.5\%$ is achieved on the test set. To demonstrate the prediction capability of the GP model, we conduct model-based forward simulation using the GP model. In this test, the GP model is used to predict the melt pool area over time for Case \#2 in Table~\ref{table_tests}. Fig.~\ref{fig_training_accuracy} shows the comparison results between the GP model prediction and the FE simulation. It is clear that the GP model predicts the FE simulation results closely.

\subsection{MPC Performance}

\begin{figure*}[htp!]
	\centering
	\subfigure[]{
		\label{fig_control_area}
		\includegraphics[width=.45\textwidth]{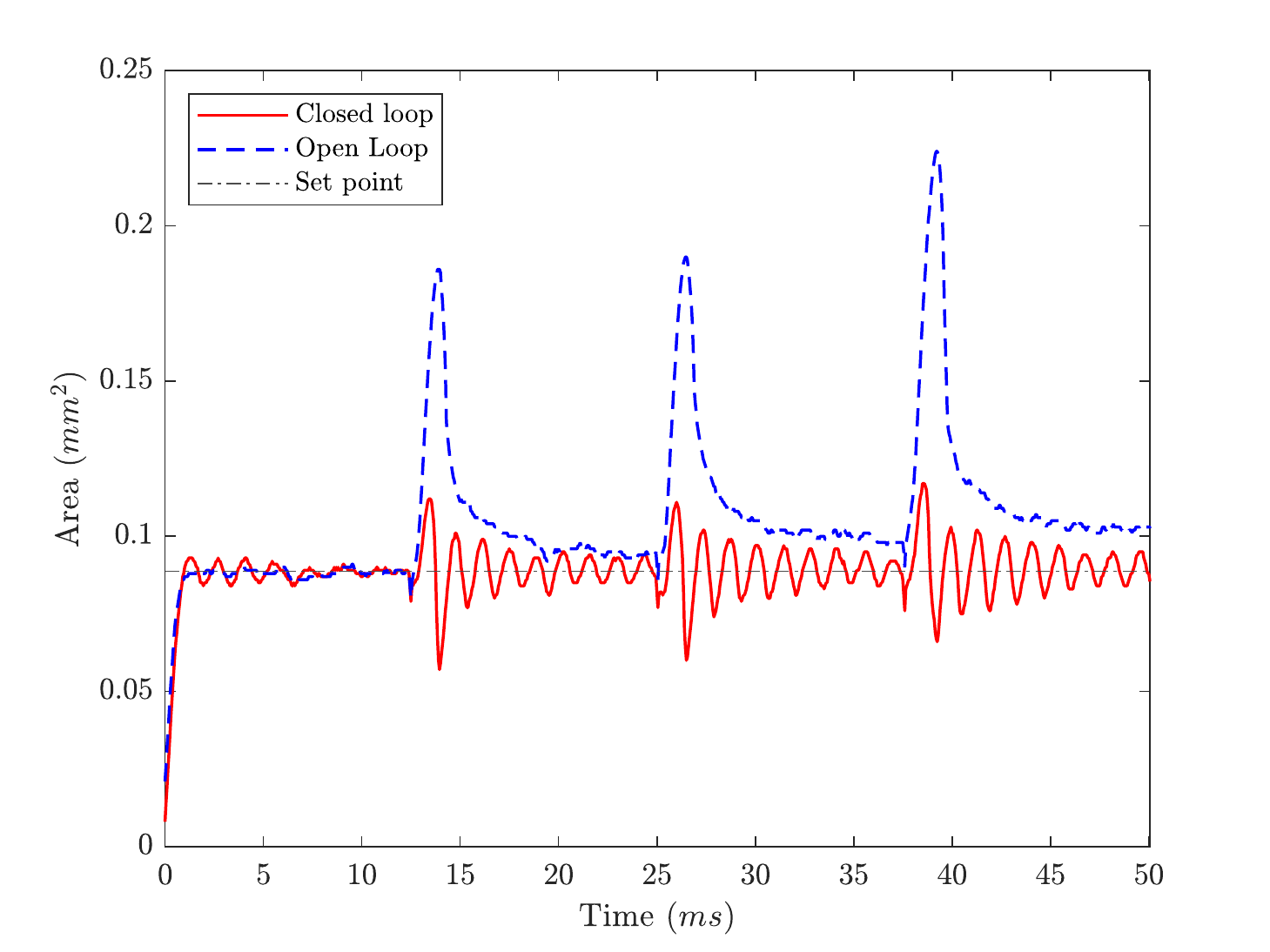}}
	\hspace{10mm}
	\subfigure[]{
		\label{fig_control_power}
		\includegraphics[width=.45\textwidth]{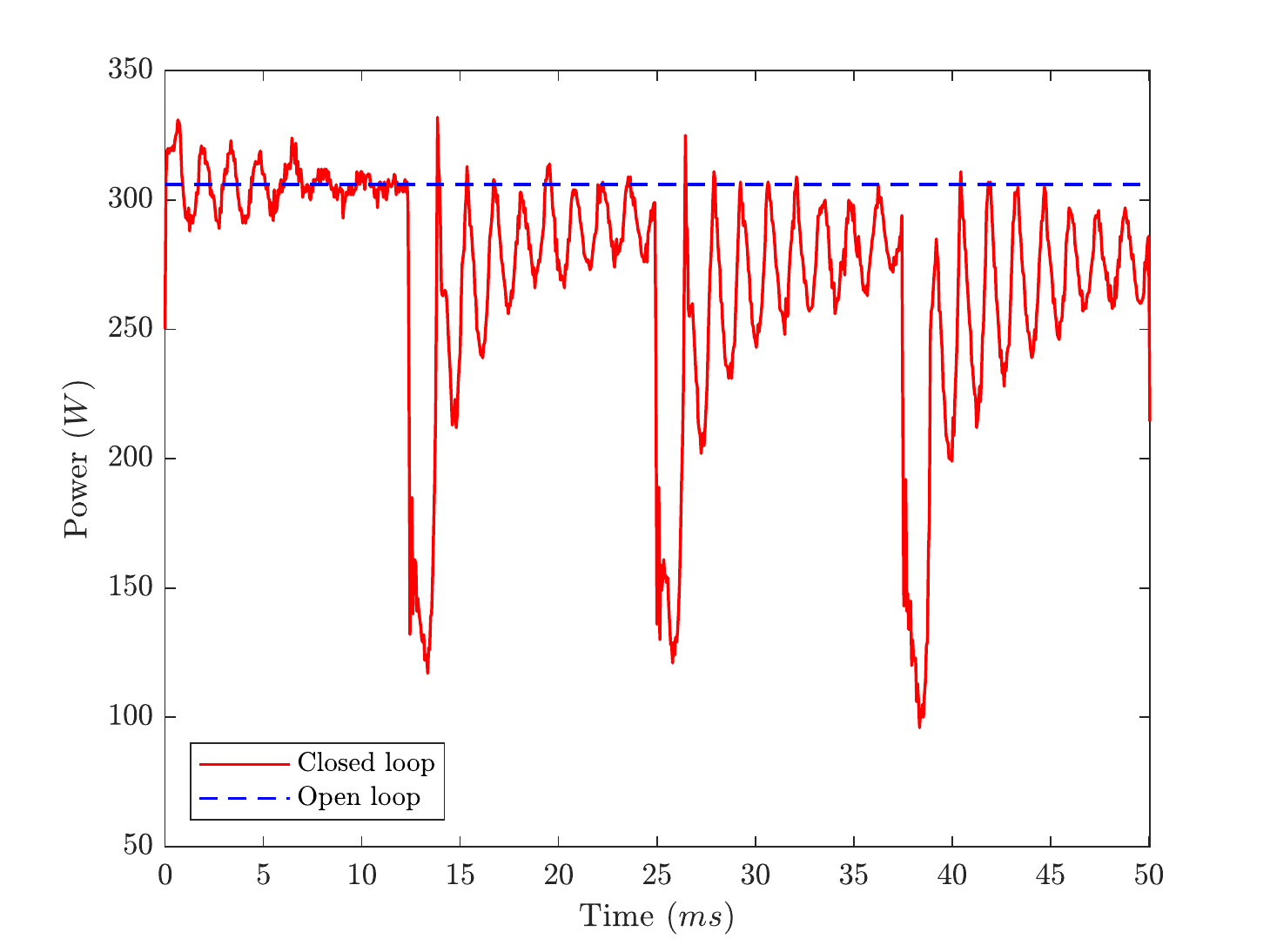}}
	\caption{Control performance comparison under MPC and constant power inputs. (a) Melt pool area $x(t)$ profiles. (b) Input laser power $p(t)$.}
	\label{test2}
\end{figure*}

Using the GP model, the MPC is implemented by only considering the laser power as the control input with constant scanning speed $v=800$ mm/s. The values of the MPC design parameters are listed in Table~\ref{table_mpc_cons}. To show the performance of the MPC controller, a $4\times 10$  mm track test is designed. The MPC controller is used to update the laser power in FE simulation to maintain constant melt pool cross sectional area over all tracks. Also, a differential plus direct feedback from the initial temperature is added to the MPC command. The challenge in multi-track SLM process is the high initial temperature after first track is completed; see Fig.~\ref{fig_test_init_temp}. The increased initial temperature leads to increased melt pool cross sectional area and therefore, feed forward control of the process is inefficient and closed-loop control is needed. Fig.~\ref{test2} shows the transition of the melt pool area and laser power profiles under the MPC design. Open-loop performance (under constant laser power input) is also included for comparison purpose. 

\renewcommand{\arraystretch}{1.3}
\setlength{\tabcolsep}{0.032in}
\begin{table}[htp!]
	\centering
	\caption{Values of the MPC design parameters}
	\label{table_mpc_cons}
	\begin{tabular}{|c|c|c|c|c|c|c|c|}
		\hline\hline 
		Parameter & $p_{l}/p_u$ (W) & $\Delta p_{l}/\Delta p_u$ (W) & $x_l/x_u$ (mm$^2$) & $Q$ & $R$ & $Q_f$ & $H$ \\ \hline
		Value & $0/350$ & $-350/350$  & $0/0.5$	 & $1$ & $0.1$ & $20$ & $20$ 		\\ \hline\hline
	\end{tabular}
\end{table}
\renewcommand{\arraystretch}{1.3}
\setlength{\tabcolsep}{0.05in}
\begin{table}[htp!]
	\centering
	\caption{Control performance}
\label{table_mpc_perf}
	\begin{tabular}{|c|c|c|c|}
		\hline\hline  
		Set point& Max overshoot & Max undershoot& Steady state $\left|e\right|_{avg}$\\ \hline
		$0.09$ & $30\%$ & $36.37\%$ & $\pm 8.16 \%$\\
		\hline\hline
	\end{tabular}
\end{table}

Despite the oscillation in the closed-loop response, it can be seen that the MPC successfully attenuates melt pool overshoot in start of new tracks, that is, at around $13$, $26$, $39$ ms in Fig.~\ref{fig_control_area}. Also, as shown in Fig.~\ref{fig_control_power}, the power level decreases significantly at the beginning of each new track. This is due to energy accumulation of multiple tracks. These results agree with these reported in~\cite{WangAM2020} that used the same material and feedback signal in their work. Furthermore, after the sudden drop at the location of each new track, the gradual increase of the laser power along each individual track is also in accordance with the experiments in~\cite{WangAM2020}. This gradual increase happens due to the exponential drop of initial temperature along each new track; see Fig.~\ref{fig_test_init_temp}. It is obvious that under the MPC closed-loop control, the variation of the melt pool cross sectional area profiles shows consistently stable results. Table~\ref{table_mpc_perf} further summarizes the closed-loop performance. Comparing with the control performance in~\cite[Table 5, five track case]{WangAM2020}, the overshoot of the response under the MPC is $10$ \% higher, the undershoot is however reduced by $40$ \%. The max overshoot and undershoot are calculated as $\frac{A_{max, min}}{A_{\textit{set point}}}$. The average steady state error is also calculated as root mean square error (RMSE) between melt pool area and the set point.
The computational burden of the control algorithm is really important in SLM process, due to high speed nature of the SLM process. The proposed algorithm takes on average $1\:ms$ to be computed. This result is achieved in python, run on a core i5\textsuperscript{\textregistered} laptop with 4 GB ram. This ensures that the algorithm can be run on real time devices with required control frequency. Note that the simulation results in the paper are done with higher sampling frequency due to limitations in the FEA time length.

\section{Conclusion and Future Work}
\label{concl}

In this paper, a data driven GP model was proposed for the SLM process. The modeling development was also used in MPC framework to regulate the melt pool cross sectional areas in multi-track SLM scanning fabrication. The main challenge in multi-track scanning scheme was the melt pool area due to the high initial temperature of substrate generated by previous scans. We used high-fidelity FE simulation that has been validated by experiments (see \cite{OlleakJAMT2020}) to generate the training data sets for the GP models. The GP model and the MPC performance were demonstrated to successfully maintain a stable, consistent melt pool cross sectional area during multi-track scanning process. Comparison with the existing feed forward SLM control was also presented. We are currently trying to implement the proposed GP model and real-time control on an SLM machine experimentally. Inclusion of other real-time sensing and laser scanning speed control is also among the ongoing research directions.   

%\begin{comment}
\section*{Acknowledgments}

The authors would like to thank Ansys\textsuperscript{\textregistered} for providing the software, and Professor Zhimin Xi of Rutgers University for his helpful discussions and comments.
%\end{comment}

\bibliography{refs}

\end{document}